\documentclass[12pt]{amsart}

\usepackage{amsfonts,amssymb,amsmath,amscd}
\usepackage{graphicx}

 \newcommand\blfootnote[1]{%
     \begingroup
     \renewcommand\thefootnote{}\footnote{#1}%
     \addtocounter{footnote}{-1}%
      \endgroup
    }

\begin{document}

\title{Robert J.~Daverman (1941 -- 2026)}

\author{Craig Guilbault}

\author{Gerard Venema}

\maketitle

\blfootnote{This is an expanded version of an article that is to appear in the \emph{Notices of the American Mathematical Society}.}
\blfootnote{Craig Guilbault is Professor of Mathematics at the University of Wisconsin-Milwaukee.  His email address is  craigg$@$uwm.edu.}
\blfootnote{Gerard Venema is Professor of Mathematics Emeritus at Calvin University.  His email address is venema$@$calvin.edu.}

Robert Jay Daverman (“Bob” to all who knew him) passed away on February 9, 2026. 

\begin{figure}[ht]
\centerline{\includegraphics[width=.85\textwidth]{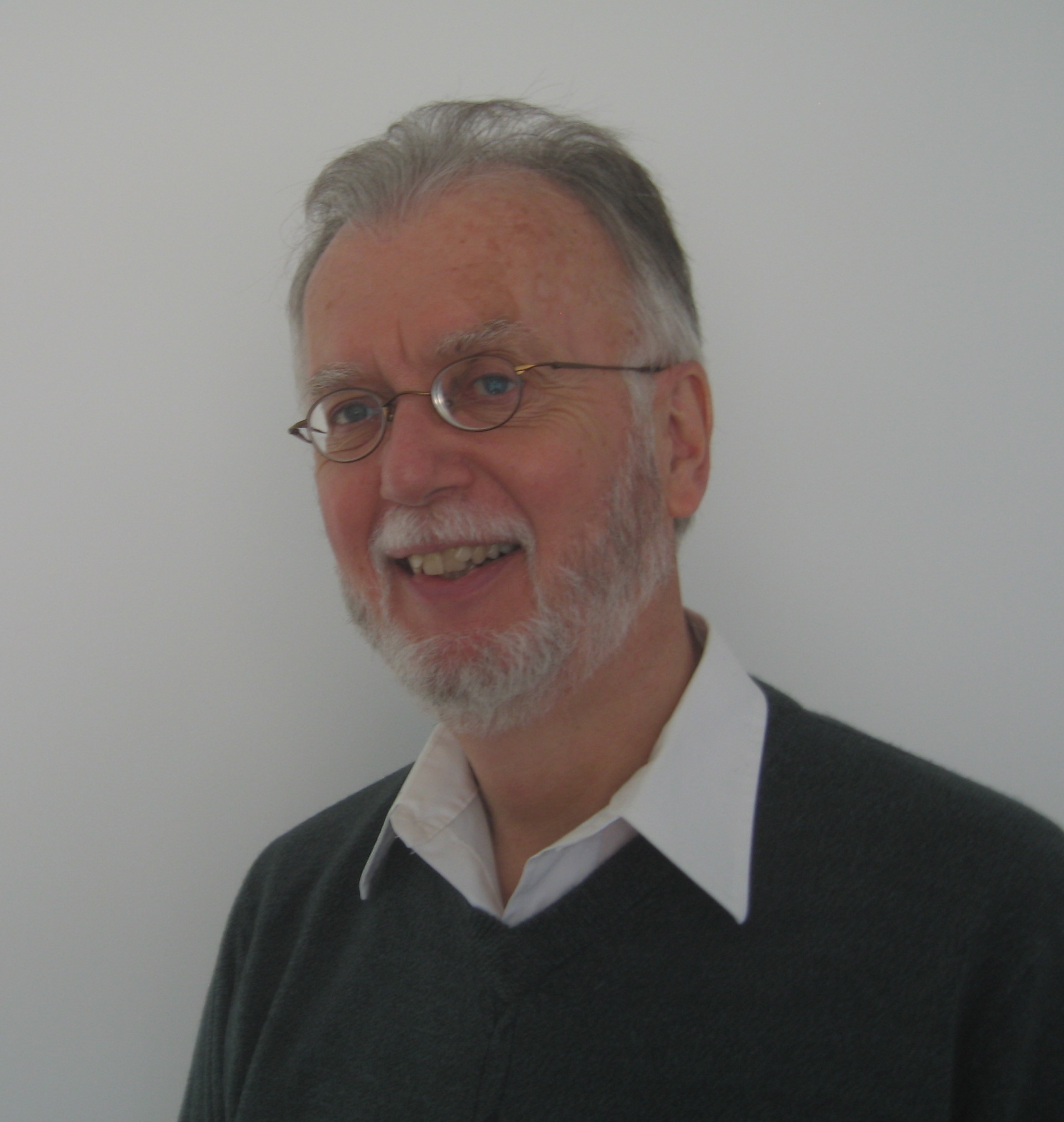}}
\caption{Bob Daverman}
\end{figure}

One of three children born to Herbert and Marian Daverman, Bob Daverman came into the world in the fall of 1941 in Grand Rapids, Michigan. That is where he grew up, met his wife and lifetime partner Lana, and earned an undergraduate degree in mathematics at Calvin College (now Calvin University) in 1963. The next stop for Bob was graduate school at the University of Wisconsin. His time there fell within the early years of an era regarded by many as a golden age in manifold topology. Under the leadership of R.H. Bing, Madison was a hotbed of activity. During just Bob’s years in Madison (1963-1967), twelve students earned PhDs under Bing’s direction---with numerous others coming before and after. Of the many successful mathematicians emerging from the “Bing school of topology,” Bob Daverman was one of the most prolific and influential.

\begin{figure}[ht]
\centerline{\includegraphics[width=.6\textwidth]{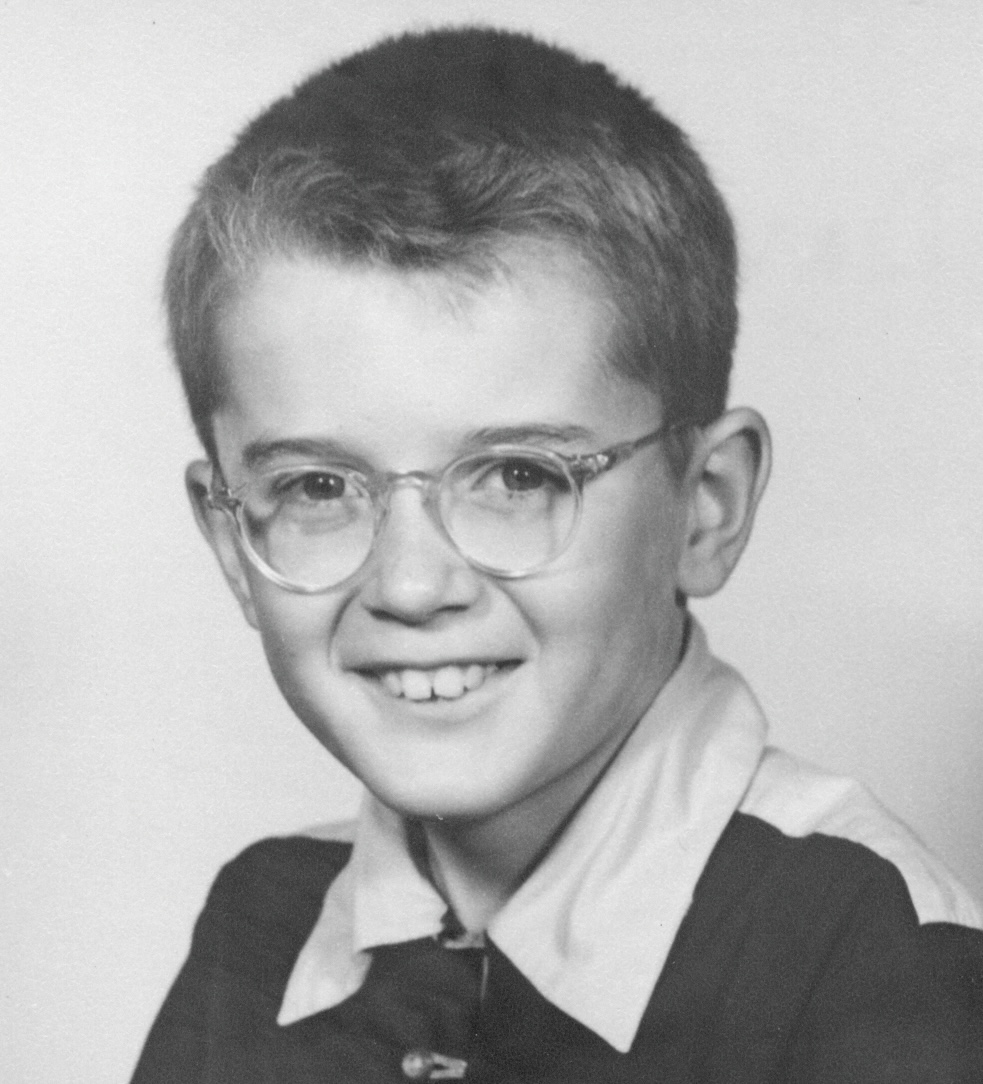}}
\caption{Young Bob}
\end{figure}

Upon graduation, Bob and Lana, with their young children Lara and Kurt, moved to Knoxville where Bob joined the faculty at the University of Tennessee---the institution where he would spend his entire academic career. It was there, through his teaching, research, mentoring, and authoring of books, that he left the largest mark on his field and on the lives of the many individuals whose paths crossed his. 

Bob'a work as a professor of mathematics was later supplemented by a parallel career of service to the American Mathematical Society where he served as Associate Secretary for the Southeastern Section from 1993 through 1999, as AMS Secretary from 1999 through 2012, and as a long-time member of the AMS Council and Executive Committee.  With his blend of professionalism and personal skills, Bob made a significant contribution to the health of the global mathematical community.

\begin{figure}[ht]
\centerline{\includegraphics[width=.9\textwidth]{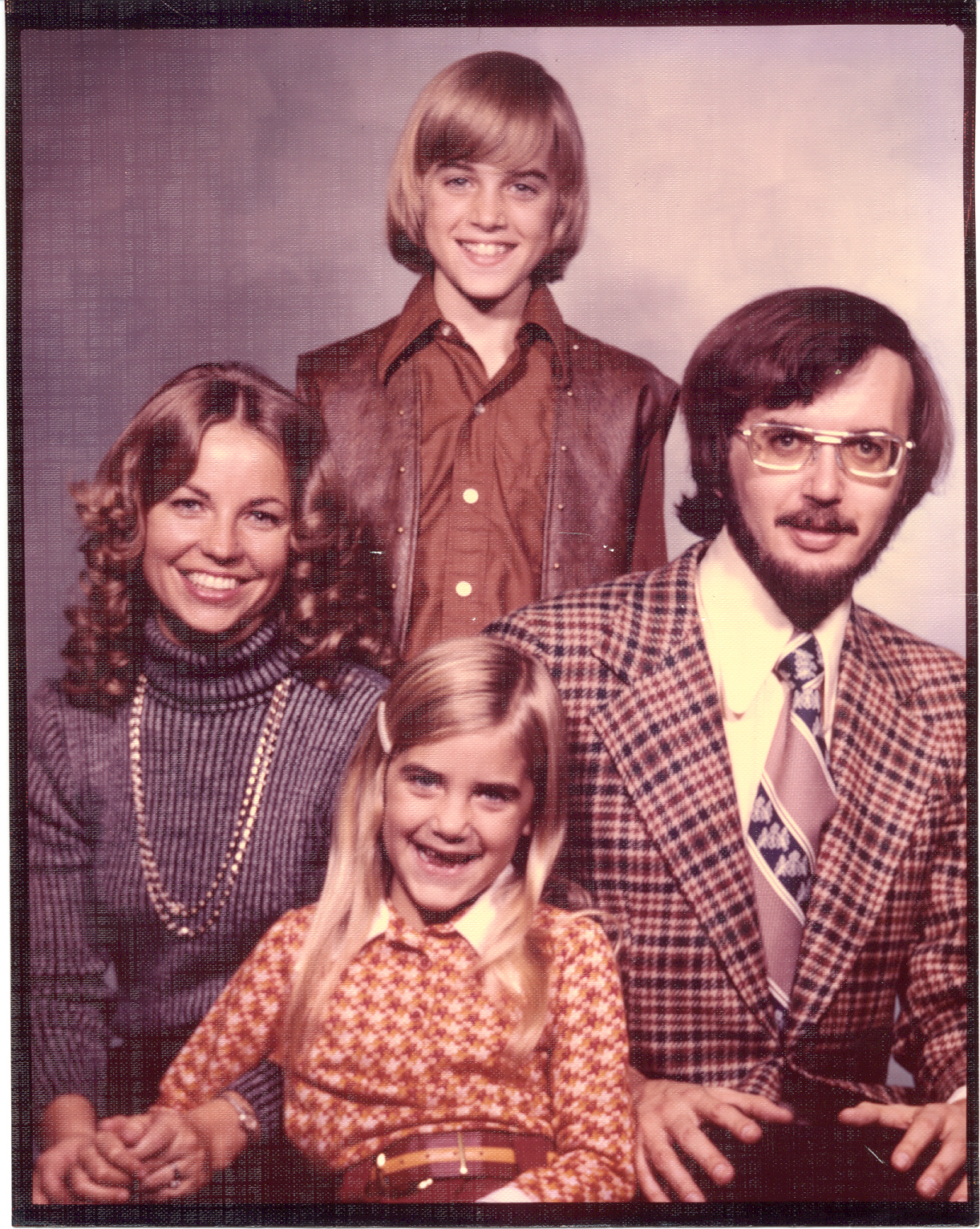}}
\caption{Bob's family: Lana, Lara, Kurt, and Bob}
\end{figure}

\section*{Mathematical Research}

\begin{figure}[ht]
\centerline{\includegraphics[width= 1.0\textwidth]{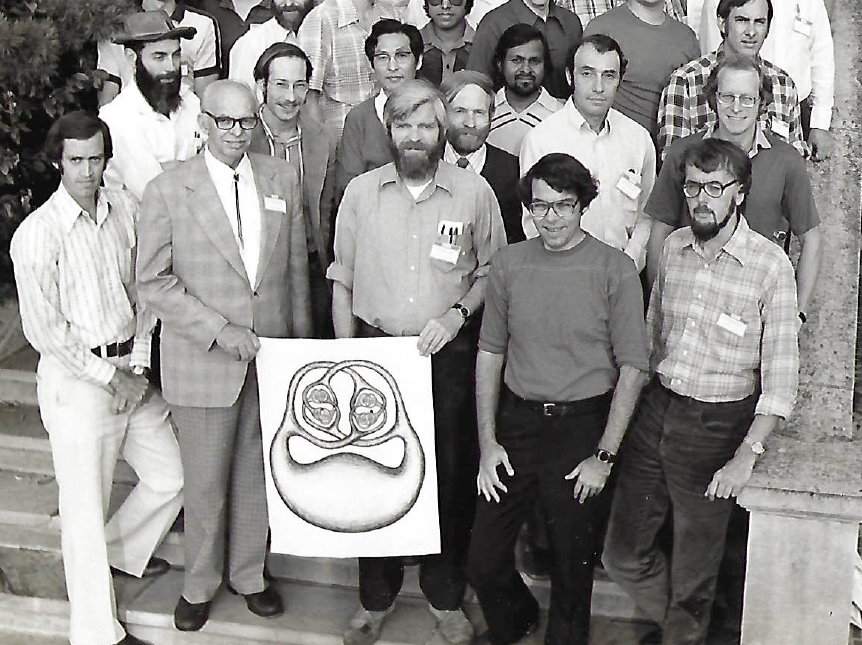}}
\caption{Bing, Edwards, and Daverman at the 1978 CBMS conference at Oklahoma State University.  Bing and Edwards are holding a drawing of the Alexander Horned Sphere.  Daverman is front right.  The others in the front row are Jim Maxwell and Dick Sher.
Larry Siebenmann is the first person in the second row.}
\label{OKphoto}
\end{figure}

When Bob Daverman began his research career as a graduate student under Bing's direction, there was a great deal of interest in wild topological embeddings.  A \emph{topological embedding} of a space $X$ into a space $Y$ is a homeomorphism of $X$ onto a subspace of $Y$.  Two embeddings $e_1,e_2:X\to Y$ are \emph{equivalent} if there is a homeomorphism $h:Y\to Y$ such that $e_2=h\circ e_1$.  An embedding $e:S^{n-1}\to S^n$ is \emph{flat} (or \emph{tame}) if it is equivalent to the standard inclusion $S^{n-1}\hookrightarrow S^n$ and is \emph{wild} otherwise.\footnote{A classic example of a wild embedding is the Alexander horned sphere---see Figure~\ref{OKphoto} where R.H.\ Bing and Bob Edwards are shown holding a drawing of the Alexander sphere.}  

\begin{figure}[ht]
\centerline{\includegraphics[width=.9\textwidth]{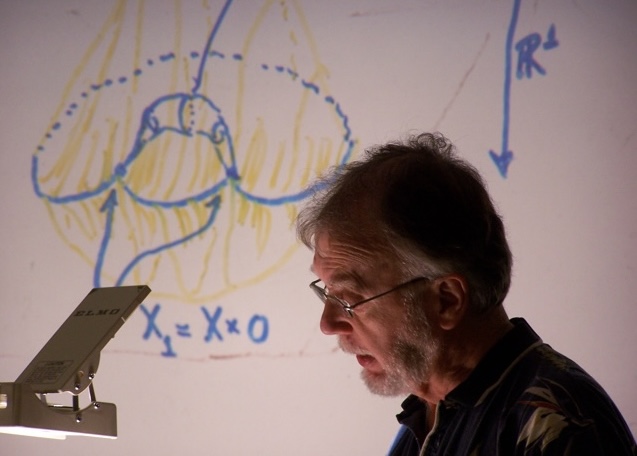}}
\caption{Bob speaking about his research}
\end{figure}

An embedding $e:N\to M$ of manifolds is \emph{locally flat} if each point $x\in e(N)$ has a closed neighborhood $U$ such that the pair $(U,U\cap e(N))$ is homeomorphic to a standard ball pair.  The Generalized Sch\"onflies Theorem\footnote{Proved by Morton Brown and, independently, by Barry Mazur and Marston Morse.} asserts that every locally flat embedding $e:S^{n-1}\to S^n$ is flat. A fundamental result of Bing is that an embedding of a surface in a 3-manifold is locally flat if the complement of the embedding is uniformly locally simply connected.

One of Daverman's lifelong interests was the study of complements of topologically embedded spheres.  The image of a topological embedding $e:S^{n-1}\to S^n$ separates $S^n$ into two connected complementary domains.  The closure of each is called a \emph{crumpled cube}.  The boundary of a crumpled cube is always an $(n-1)$-sphere and two crumpled cubes can be sewn together by identifying points that correspond under some homeomorphism of the boundaries.  In his early days at Tennessee, Bob proved several basic results about sewings of 3-dimensional crumpled cubes, some in collaboration with his colleague Bill Eaton.  For example, given any two crumpled cubes, there is a sewing that yields the 3-sphere.   And there is an uncountable family of topologically distinct crumpled cubes that are universal in the sense that they can be sewn to any other crumpled cube in a way that yields the 3-sphere.

\begin{figure}[h]
\centerline{\includegraphics[width=\textwidth]{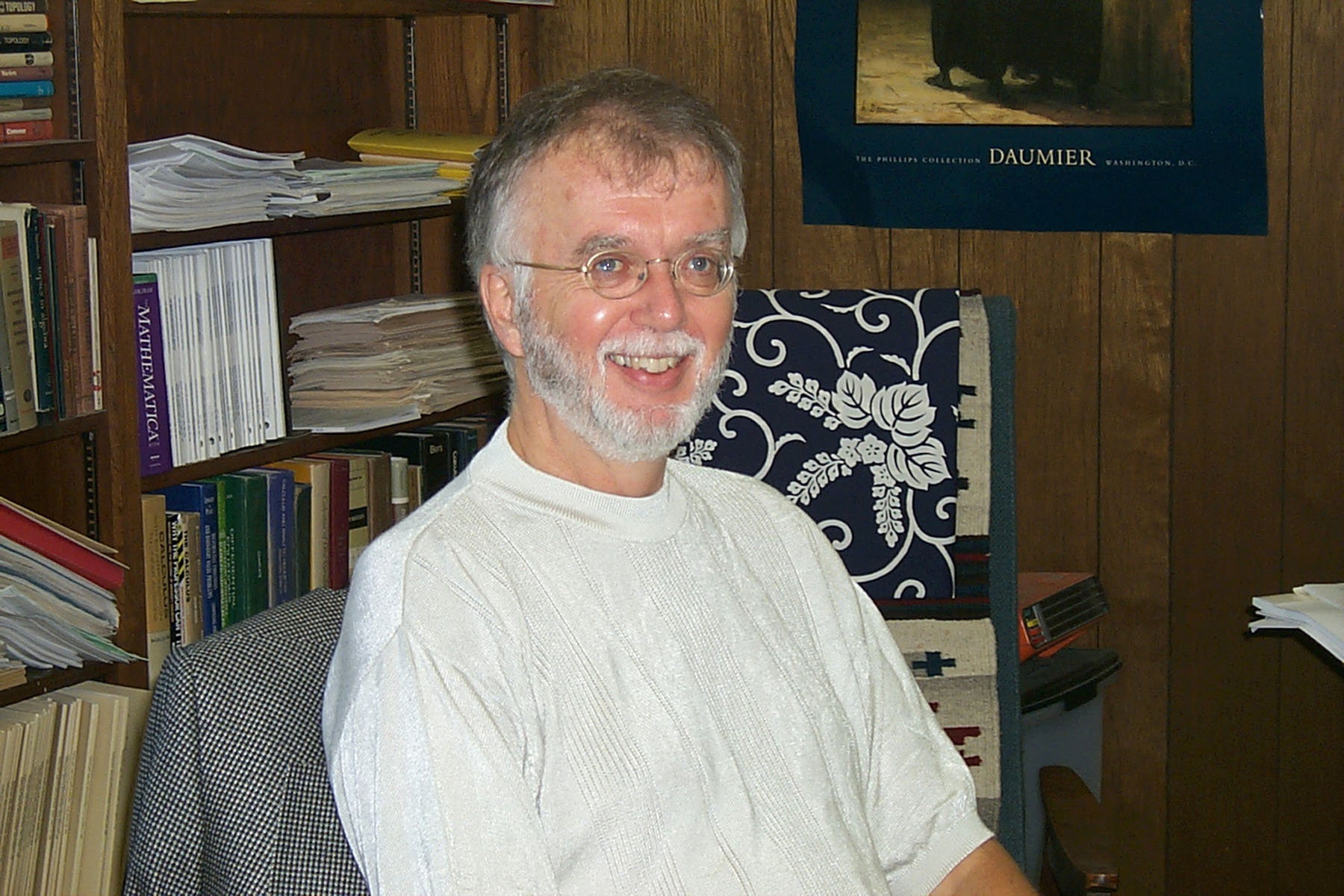}}
\caption{Bob in his office at the University of Tennessee}
\end{figure}

Around 1972, Daverman began to work on extending the theory of tame-versus-wild embeddings to higher dimensions.  One notable result from that time is his proof that an embedding of an $(n-1)$-dimensional manifold into an $n$-manifold ($n\ge5$) is locally flat provided the complement of the embedded manifold is uniformly locally simply connected, extending Bing's result to higher dimensions.\footnote{A.\ V.\ \v Cernavski also contributed to this theorem.}  Thus wildness is detected by the first homotopy group, regardless of the ambient dimension.\footnote{The theorem was later extended to the missing dimension 4 by Freedman and Quinn.}  Other significant results include a slicing theorem for embeddings of spheres and a surprising example that shows the necessity of ``twice-tameness'' in Rob Kirby's theorem regarding embeddings that are tame modulo a Cantor set.  Daverman became the authority on high-dimensional crumpled cubes.  He characterized what types of crumpled cubes are theoretically possible and constructed examples illustrating the various types. In collaboration with Jim Cannon, he proved a noteworthy theorem about sewings of high-dimensional crumpled cubes known as the Mismatch Theorem.

During the 1970s Daverman increasingly turned his attention to decomposition theory.  A \emph{decomposition} of a metric space $X$ is simply a partition of $X$ into compact subsets.  Associated with a decomposition $G$ there is a \emph{decomposition space} $X/G$, whose points are the elements of $G$ and which is endued with the quotient topology.  It is always assumed that the decomposition is \emph{upper semicontinuous} (usc), a technical condition which ensures that the decomposition space is metrizable.  In decomposition theory the properties of a surjective map $f:X\to Y$ and the image space $Y$ are investigated by studying the geometric properties of the elements of the decomposition $G=\{f^{-1}(y) \mid y\in Y\}$.  

One basic problem in decomposition theory is to determine conditions under which the projection map $\pi:X\to X/G$ can be approximated by homeomorphisms.  The \emph{Bing Shrinking Criterion} gives a simple, elegant answer:  $\pi$ can be approximated by homeomorphisms if and only if it is possible to simultaneously shrink all the elements of $G$ to small size in $X$ with a homeomorphism $h:X\to X$ for which $\pi\circ h$ is close to $\pi$.  

In case the domain space is a manifold, a necessary condition for the decomposition to be shrinkable is that the decomposition elements be \emph{cell-like}, which means that each decomposition element is contractible in a small neighborhood of itself.  A major question at the time was what other conditions are needed to ensure that a decomposition of a manifold is shrinkable.  The answer was provided by Bob Edwards,\footnote{Edwards writes. ``I owe a great deal to Bob Daverman for conversations we had regarding developments in geometric topology.  Especially important were those involving Bing-type decomposition theory, that took place when we were both visiting the University of Utah in 1974, and which were critical to my work in the years that followed.''} who proved the following breakthrough result in 1977: If $G$ is a usc decomposition of an $n$-manifold $M$, $n\ge5$, into cell-like sets, then the decomposition map $\pi:M\to M/G$ can be approximated by homeomorphisms if and only if $M/G$ is finite dimensional and satisfies the disjoint disks property.  A space $X$ satisfies the \emph{disjoint disks property}  (DDP) if any two maps of the 2-cell into $X$ can be approximated by maps with disjoint images.\footnote{The DDP was first identified by Jim Cannon, who had proved a significant special case of Edwards' Theorem earlier in 1977.}

\begin{figure}[ht]
\centerline{\includegraphics[width= 1.35\textwidth]{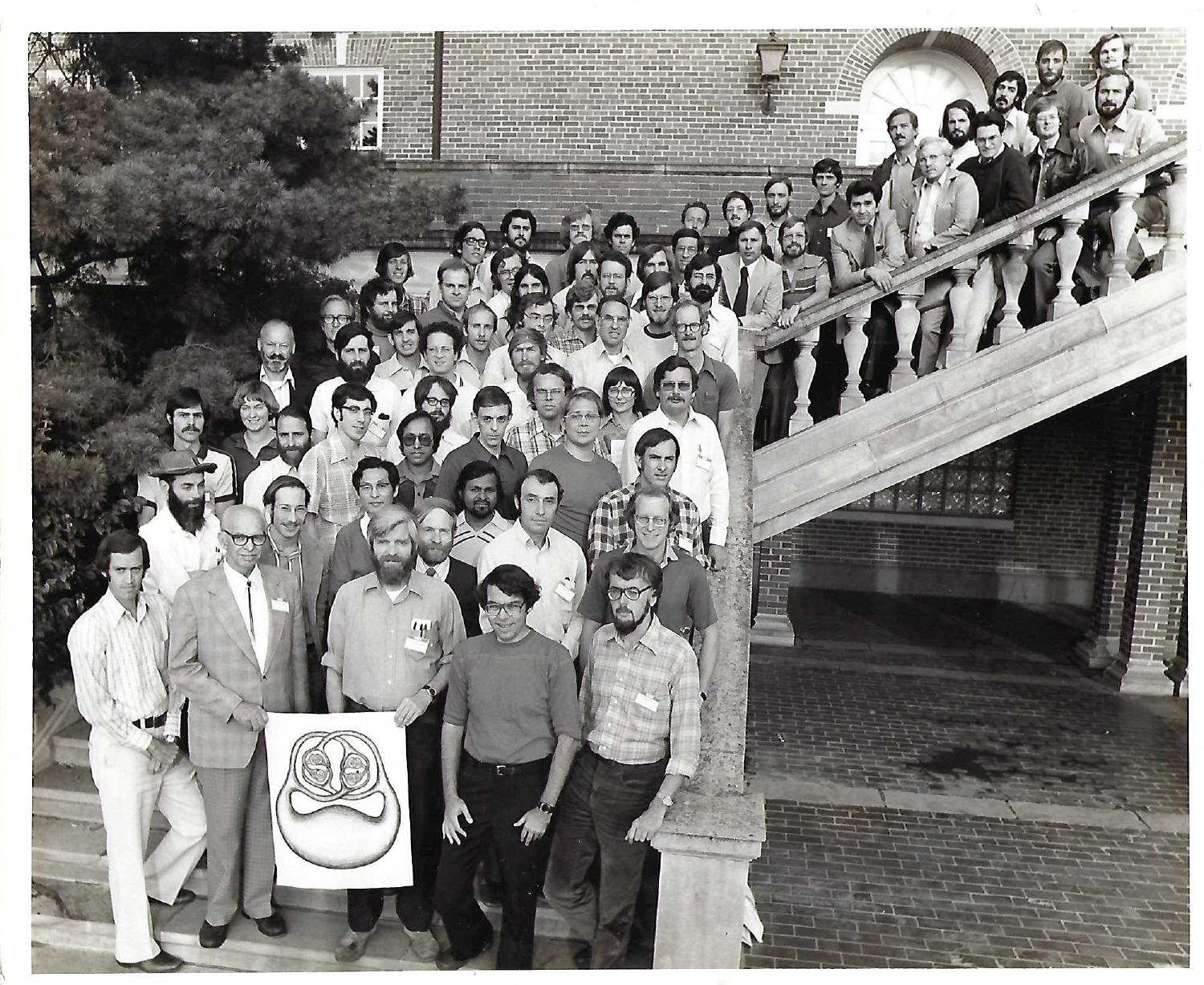}}
\caption{All the participants in the 1978 CBMS conference at Oklahoma State University}
\end{figure}

Over the next few years Daverman produced dozens of new results about decompositions of manifolds, many based in some way on Edwards' work.  One such result is an example of a non-shrinkable decomposition of $S^n$ ($n\ge5$) whose non-degenerate elements form a null sequence of cellular sets, which is a high-dimensional analogue of Bing's famous dog bone space.\footnote{A compact subset of an $n$-manifold is \emph{cellular} if it has small neighborhoods that are $n$-cells.  A set can be cell-like (an intrinsic property) without being cellular (a property of the embedding).}  Daverman proved that stabilizing a decomposition of a manifold by crossing with $\mathbb{R}^2$ always makes the decomposition shrinkable.  He also proved that in many cases crossing with $\mathbb{R}^1$ is enough to make the decomposition shrinkable, but whether that is true in general remains an open problem.

Daveman's work in decomposition theory cannot be separated from his work on embeddings.  He quotes Jim Cannon as saying that there is ``an intimate connection between decomposition space theory and taming theory, a connection enriching and unifying both areas.'' Two early examples that point to the interrelatedness are Bing's proof in the early 1950s that the sum of two solid Alexander horned spheres\footnote{The solid Alexander horned sphere is the crumpled cube determined by the non-simply connected complementary domain of the Alexander sphere.} is $S^3$ and Brown's proof of the Generalized Sch\"onflies Theorem.  Exploiting this connection, Daverman constructed decompositions that produced  surprising examples of embeddings of the Cantor set and cells of various dimensions. In joint work with John Walsh he described an example of a ``ghastly generalized manifold'' which contains no 2-disks whatsoever.  Daverman and Cannon constructed a decomposition that yields an example of a flow in which every flow line is wild at every point.

In the summer of 1977 Edwards circulated a hand-written write-up of the proof of his Cell-like Approximation Theorem and he also gave a complete exposition of his theorem and its proof at the 1978 CBMS conference held at Oklahoma State University.  Edwards never submitted his manuscript for publication,\footnote{Recently Ric Ancel and the authors of this article have edited the manuscript and it will appear in \emph{Celebratio Mathematica}.} but Daverman included the proof in his  book \emph{Decompositions of Manifolds} \cite{DecompBook}, which is a comprehensive exposition of the theory of decompositions and which became the go-to reference for those wishing to understand the two critical shrinking arguments in Freedman's proof of the 4-dimensional Poincar\'e Conjecture.\footnote{The proof of the Cell-like Approximation Theorem requires special care in dimension~5.  Daverman and his student Denise Halverson published the details in 2007.} 

\begin{figure}[h]
\centerline{\includegraphics[width= .85\textwidth]{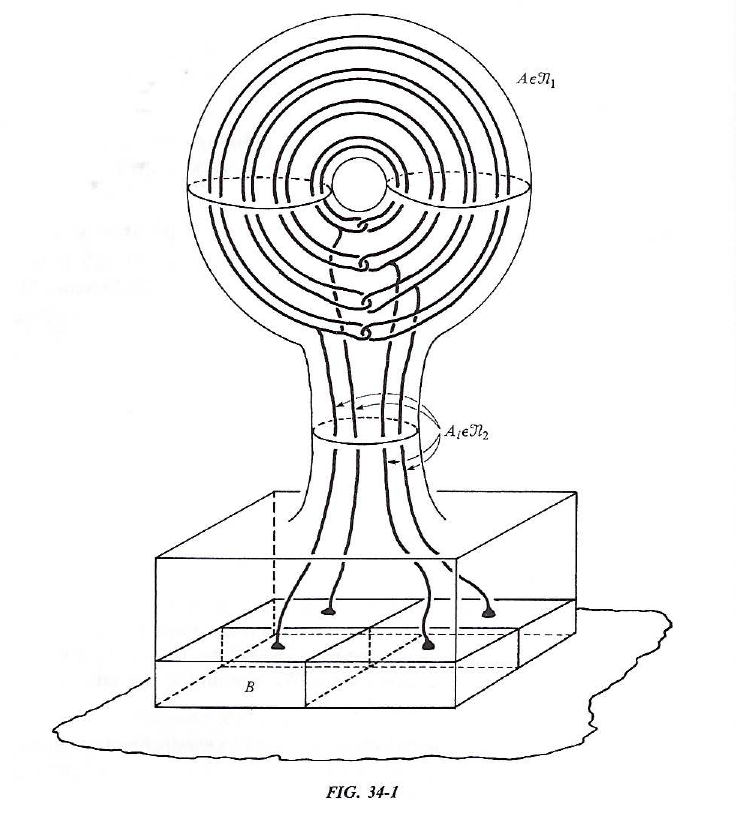}}
\caption{Fig.~34.1 from \cite{DecompBook}.  Bob grew up in a family of architects and had a particular talent for technical illustration}
\end{figure} 

Around 1990 Daverman introduced and defined a class of manifolds he called fibrators.  A closed manifold $N$ is a \emph{fibrator} if for every usc decomposition of a manifold such that each element of the decomposition has the homotopy type of $N$, the projection map is an approximate fibration.  Over the next several decades he published a series of papers in which he proved that various manifolds are fibrators.  For example, any simply connected manifold is a fibrator, as is complex projective space, while real projective space is not a fibrator. With his students Young Ho Im and Yongkuk Kim he studied fibrator properties that are preserved under products and connected sums.  His final student, Violeta Vasilevska, wrote a dissertation on fibrators in the piecewise linear category.

A \emph{lamination} of an $n$-dimensional manifold is a usc decomposition whose elements are $(n-1)$-manifolds.  Bob proved that the decomposition space of a lamination is always a 1-manifold, even if the submanifolds in the decomposition are wildly embedded.  He continued his study of laminations in a series of joint papers with Fred Tinsley in which they eventually found (crumpled) laminations of the cobordisms in Quillen's plus construction.

Daverman coauthored, with Gerard Venema, a second book, titled \emph{Embeddings in Manifolds} \cite{EmbBook}. The book focuses on embeddings of compact polyhedra in piecewise linear manifolds, particularly on the questions of which embeddings are tame and which are equivalent.  (In this context, an embedding is \emph{tame} if it is equivalent to a piecewise linear embedding.)  The book contains descriptions and illustrations of many of the classical examples of topological embeddings.  The theorems are organized by codimension, which is the difference between the dimension of the embedded polyhedron and that of the manifold in which it is embedded.  In each codimension appropriate conditions on the fundamental groups play an important role.

Daverman maintained his interest in crumpled cubes throughout his research career.  Four of his last seven research papers address questions about crumpled cubes.  In his last paper on the subject, he and Shijie Gu propose a hierarchy for crumpled cubes.  The wildness of two crumpled cubes can be compared if there is a usc decomposition that projects one onto the other. 

A prolific researcher, Daverman published well over 100 research articles and two research monographs (\cite{DecompBook} and \cite{EmbBook}).  He contributed to research in geometric topology by editing \emph{The Collected Papers of R.H.\ Bing} \cite{BingPapers} and \emph{The Handbook of Geometric Topology} \cite{Handbook}).  Two of his special gifts were asking good questions that stimulated and advanced research in his field (for instance, his problem list in \cite{Problems}) and producing creative and illuminating examples.

\section*{Mentoring}

To former students, twelve of whom earned PhDs under his direction, Daverman was a towering presence at the University of Tennessee, Knoxville (UTK). As a professor, his teaching stood out. True to his roots in the Bing school of topology, he employed the ``Moore method'' in introductory graduate topology classes. In place of a text, there was a terse set of notes containing definitions and propositions but no proofs. The game, for students, was to figure out the proofs and present them to the class. Daverman played the dual role of coach and referee. For UTK graduates of various specialties, that course has been frequently recalled as a formative mathematical experience. For those who specialized in topology, it was just the beginning. In advanced topics classes, Bob gave crisp, well organized lectures bursting with colored chalk. Embedding theory and decomposition theory were specialties but, at the request of students, he would offer a class on just about any relevant topic involving manifolds—a task that surely involved enormous time commitments. 

Seminars were another area where Bob showed exceptional generosity toward his students. In a department where they did not count toward faculty workloads, seminars played a central role in the education of topology students. Along with the standard research seminar in topology, there were student topology seminars. These typically involved a handful of graduate students, meeting twice weekly and working through a foundational text (for example, \emph{Dimension Theory} by Hurewicz and Wallman or \emph{Piecewise-linear Topology} by Rourke and Sanderson). Bob attended and provided guidance throughout.

Daverman’s enthusiastic love for the greater mathematical community was always on full display, and he strove to bring aspiring young mathematicians into the fold.  A favorite way of introducing them to the culture was through the telling of stories involving mathematicians and mathematical events. Quirky personalities were to be enjoyed and accepted; exceptional talent was to be admired and appreciated. Another avenue for initiation into the broader community involved conferences. At a time when funding priorities were not as student friendly as they are today, Bob was ahead of his time. His efforts at inclusion often took the form of long road trips. A drive from Knoxville to Tallahassee, Florida or Lafayette, Louisiana was not out of the question. (The Spring Topology Conference was a favorite event.) By choosing a university van over an airline ticket, Bob could afford to bring the students. Those students were often among a select few in attendance—a circumstance that facilitated quick immersions into the established crowd. Afterwards, the trip home was ideal for retelling old stories (now with faces to associate with the characters) and for crafting some new ones.  

Bob’s mentoring legacy extended well beyond his own students. His outgoing personality resulted in a large network of friends in the mathematical community but, to a significant subset of that network, the connection went beyond ordinary friendship. Bob was a major force in the careers and lives of many. 

\begin{figure}[t]
\centerline{\includegraphics[width= \textwidth]{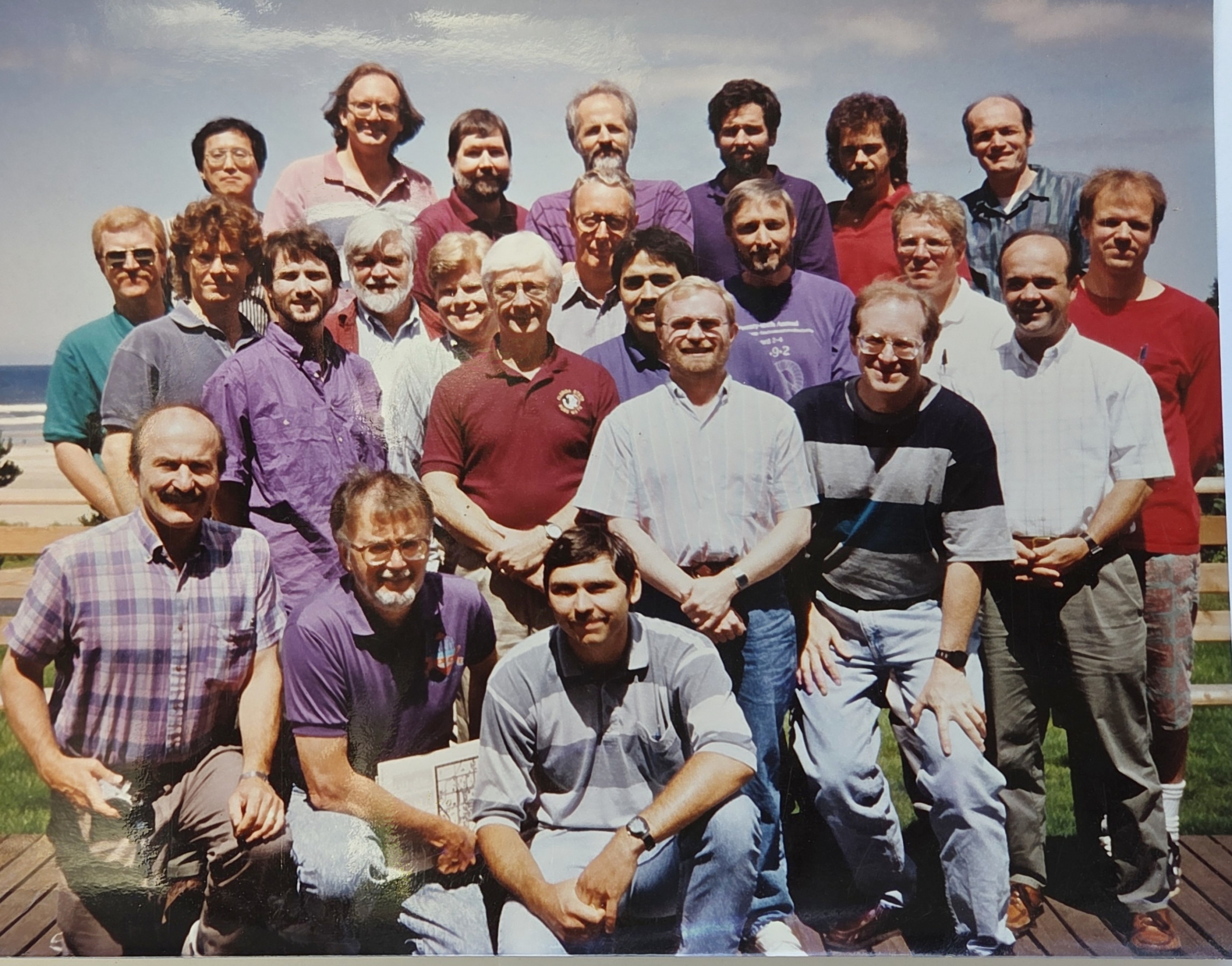}}
\caption{Participants in the 1993 Workshop in Geometric Topology, which was held at Oregon State University}
\end{figure}

A lasting symbol of his impact can be observed each June at the Workshop in Geometric Topology (WGT), which recently met for the forty-third time. At its inception in 1984, it was a small informal gathering of geographically scattered young topologists looking to share ideas and support each other in an era before easy electronic communication. Besides topology, that original group of workshop organizers: Dennis Garity; Jim Henderson; Terry Lay; Fred Tinsley and David Wright, had another thing in common—all viewed Bob as a mentor. At crucial times in the early stages of their careers, each had spent significant time at UTK, mostly arranged by Daverman, with one-on-one mathematical interaction as a primary goal. For some, those interactions were career-defining. It was, therefore, no coincidence when Bob was chosen as the first Principal Speaker at the Second Annual WGT, hosted by Colorado College. There, a tradition was established. The workshop grew but, for many, the highlight was constant: an annual opportunity for renewed interaction with their mentor. His high-energy enthusiasm for talking, encouraging, and doing mathematics were a defining feature of countless workshops. And through those workshops Bob's mentoring influence expanded.\footnote{Bob's mentoring legacy did not go unnoticed by leaders in the field. Jim Cannon observed that,``Bob served as finishing mentor for a number of my graduate students. He was very effective in helping young mathematicians make the transition from new PhD to productive mathematical citizen.''}     From its early days, a primary goal of the WGT, written into numerous successful NSF proposals, was to help mathematicians working in less-than-ideal environments—geographic isolation, high workload, or lack of financial resources—maintain active research programs. That mission blended perfectly with one of Bob’s core beliefs: that every mathematician, properly motivated, is capable of making a genuine contribution to the field. When a young mathematician showed signs of struggle or was veering off in an unproductive direction, Bob patiently took the time to offer key insights or to gently steer a course correction.  Later, when progress had been made, he was quick to lead the cheers.

\begin{figure}[t]
\centerline{\includegraphics[width= \textwidth]{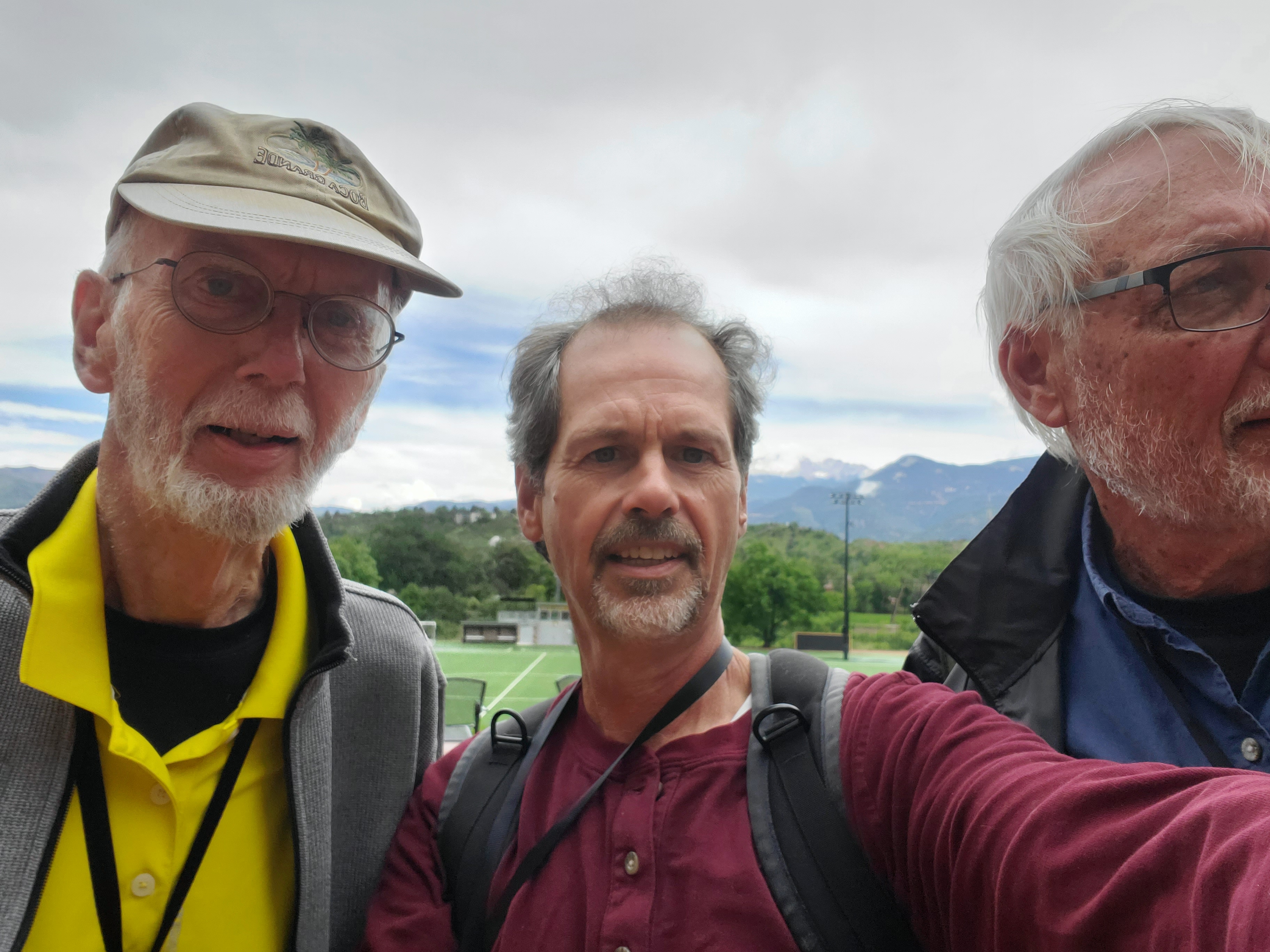}}
\caption{Daverman, Guilbault, and Venema at the 2024 Workshop in Geometric Topology, which was held at Colorado College.  It was the last conference Bob was able to attend}
\end{figure}

In subsequent years, additional Daverman mentees—his own PhD student Craig Guilbault and (eventual) co-author and fellow Calvin College alumnus Gerard Venema—were added to the WGT organizing team, further strengthening the Daverman bond. In 2002, at his alma mater in his hometown of Grand Rapids, with his family in attendance, the WGT organizers expressed their appreciation with a 60th birthday celebration added to that year’s workshop. Another commemorative event—this one in recognition of Bob's retirement—was built around a Special Session in Geometric Topology (organized by Guilbault and Steve Ferry) at the 2014 AMS Spring Sectional meeting in Knoxville. More recently, in remembrance of his passing, the 2026 WGT included special talks and a dinner where happy and poignant memories were shared.

\bibliographystyle{amsplain}
\bibliography{DavermanRefs}

\vskip1truein

\section*{Daverman's Service to the AMS}

\begin{center}
by \textsc{John Ewing}\footnote{John Ewing is President Emeritus of Math for America and former Executive Director of the AMS.  His email address is john$@$ewings.org.}
\end{center}

\medskip

Bob Daverman served the AMS for nearly twenty years, in two separate but related positions, from 1993 through 2012. In this regard, he is among a select few mathematicians who served the society for two decades---Frank Nelson Cole (1896--1920), RGD Richardson (1921--1940), and Everett Pitcher (1967--1988). Among these, Bob was singular for his calm demeanor and unfaltering diplomacy. 

Bob began his service as associate secretary in 1993. There are associate secretaries, one for each section, who concurrently oversee the various meetings of the AMS—sectional (their own), national (rotating), and international (occasional). It’s a job requiring lots of attention to detail … and lots of people skills. He was sensational in both. Well organized, decisive, precise, and gently charming, even when telling someone No (which is perhaps the most important, and the hardest, task one does when organizing a meeting). The Southeastern section was a model of efficiency.

In 1999 the position of secretary came open, and Bob applied. Because of his previous five years as associate secretary he was a shoo-in. 

To a mathematician, the title of ``secretary'' may sound not very prestigious. But that’s a mistake. When the AMS was founded in 1888, it had only two officers—president and secretary. For its first 60 years, the Society was run by the secretary. The AMS president was its titular head, but the secretary  ran the Society, later joined by a treasurer. Great luminaries held the position of secretary—Thomas Fiske, Frank Nelson Cole, RGD Richardson, and JR Kline. Unlike the president, who served two-year terms, secretaries stayed around for many years. They shaped the Society.

In 1949, the AMS created the position of executive director, who in turn hired more and more staff to support a growing society, which began to do more and more things. But the secretary remained involved in almost every aspect of the Society. The secretary, along with the ``secretariat'' (four associate secretaries), created and oversaw all meetings. The secretary, with the president, guided and interpreted the Society’s scientific policies. The secretary guided the president in appointing hundreds of mathematicians to hundreds of committees every year. All this, with a single assistant. The executive director was in charge of personnel, finances, and publishing. All the other ``stuff'' fell in between. The division of responsibilities was shadowed by the bipartite governance structure, with the Council agenda controlled by the secretary, the Board agenda by the executive director, and the president chairing both.

\begin{figure}[ht]
\centerline{\includegraphics[width=\textwidth]{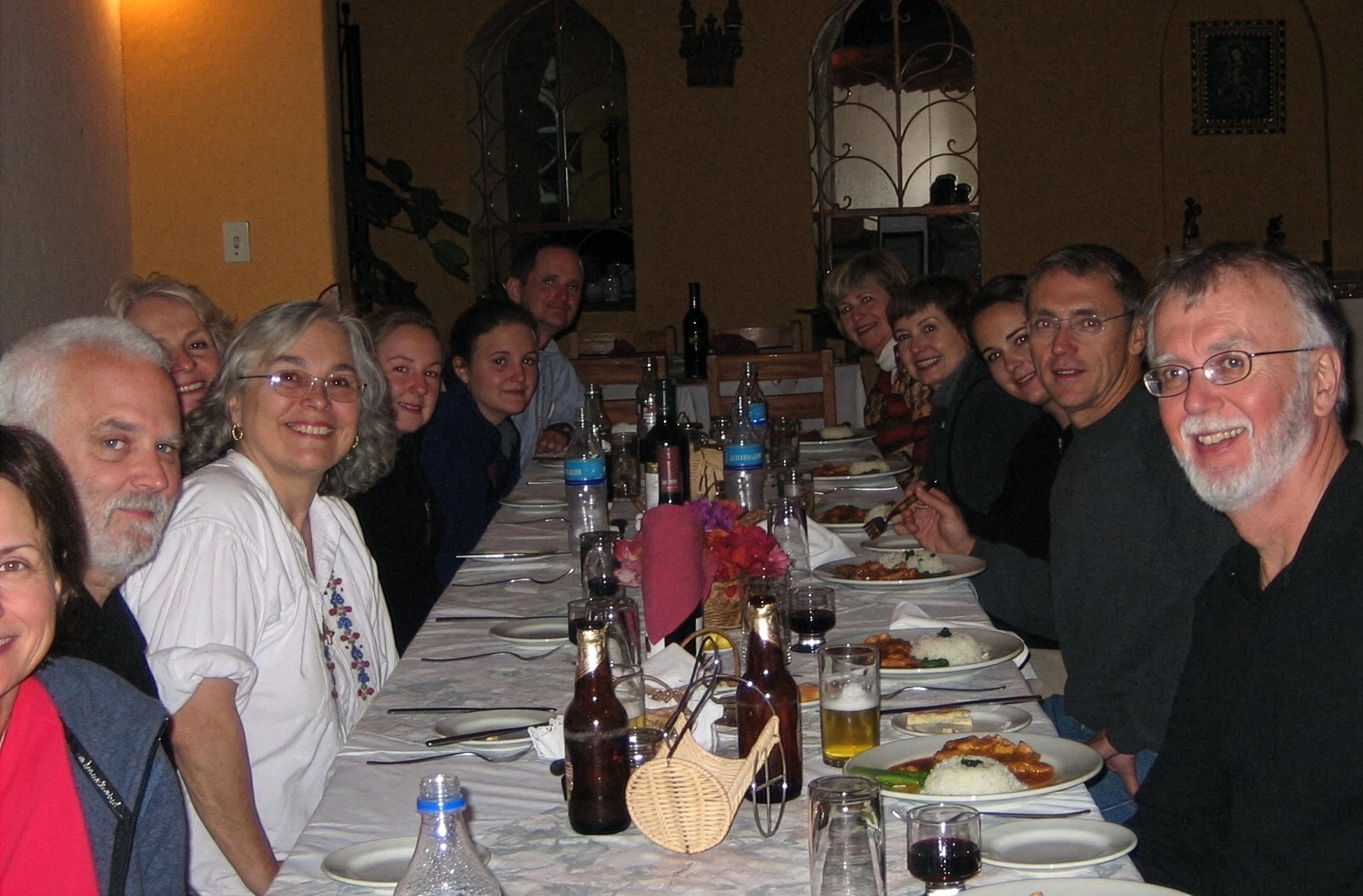}}
\caption{John and Bob at a banquet in Peru in 2007.  The photo was taken while the two of them, along with their spouses, were participating in a volunteer project that provided eye exams and glasses to thousands of rural Peruvians in the high Andes.}
\end{figure}

As the AMS continued to expand its scope during the second half of the 20th century, this arrangement grew more complex and was fraying. Able secretaries with strong commitments to the AMS worked with able presidents with strong commitments, and they worked with able executive directors with equally strong commitments. Responsibilities overlapped. Decision making grew complicated. Compromise was often required. 

I was the executive director when Bob took up his position as secretary four years into my term. He came along at a crucial time in the life of the AMS. 

Bob was an accomplished administrator, running the secretariat with exceptional efficiency and guiding the Council with care. But he was also a superb diplomat and a natural collaborator. Because he was associate secretary before his appointment, he was aware of the inherent friction in an organization with bicameral governance and tripartite leadership.  Bob set out to make things work. He built on other efforts, including those of his predecessor, but he was the most accomplished at “making it work,” no matter what the obstacles. He set about the task intentionally and enthusiastically. And he left a legacy for his successors.

His gentle demeanor worked well in other parts of his job. We served together on various committees and panels in Washington, and Bob was frequently the designated mediator. When dealing with other societies, foreign or domestic, he was the calm voice of reason. When Council members occasionally rose in anger about some issue, Bob would talk them down. (Well, usually.) His calm powers of persuasion were his hallmark, and his tenure set a precedent for the future of the AMS, right up to the present. 

Bob loved being a mathematician. He loved being associate secretary and secretary. He valued the Society’s traditions and its history. He was dedicated to his work. Over his nearly twenty years of service, he changed the AMS and left it better. For any organization, that’s the best kind of service one can provide.

\vskip1truein

\section*{Memories of Bob Daverman}
\begin{center}
by \textsc{Du\v{s}an D. Repov\v{s}}\footnote{Du\v{s}an Repov\v{s} is a professor of geometry and topology at the University of Ljubljana, Slovenia. His email address is dusan.repovs$@$guest.arnes.si.}
\end{center}

\medskip
Many mathematicians first impress us through their mathematical results. Bob Daverman certainly did that, but what I remember most clearly is the sheer pleasure he took in mathematics and the pleasure of being in his company. Within a few hours of meeting him, I understood why people were genuinely happy whenever Bob arrived.

I first met him in the late 1970s, when I was a graduate student at Florida State University and he came to Tallahassee to lecture on decompositions of manifolds. Around Bob, mathematics never remained still for long. A question produced an example, the example suggested a conjecture, and the conjecture soon met a counterexample or a new construction. All of this came with warmth, curiosity, and laughter. He gave the impression not of displaying mathematics, but of inviting everyone nearby to join in it.

After I graduated, Bob invited me to give my first talk at the University of Tennessee. He met me at the airport in his small sports car and drove me to his house, where I stayed during the visit. Our mathematical conversation began almost as soon as I got into the car and continued without interruption all the way to Knoxville. The next morning I discovered another memorable feature of life with Bob: an enormous coffee mug, surely the largest I had ever seen. Somehow he could drive at full speed, discuss topology, and drink from the mug - without spilling a drop.

Bob visited Ljubljana several times between 1986 and 1994, once as a Fulbright Fellow. His visits were very productive for the Slovenian topologists, but they also transformed the atmosphere around us. Whether we were in an office, walking through the city, sitting at dinner, or waiting for transportation, topology soon reappeared. Bob seemed to carry with him an inexhaustible supply of problems, examples, conjectures, counterexamples, and geometric pictures.

One evening, after dinner at our home, the rest of us were ready to relax. Bob was not. He went into my study and began drawing on the small blackboard. Before long, every available inch was covered with sketches of strange decomposition spaces. The blackboard remained untouched for several weeks. My children insisted that it stay exactly as Bob had left it and proudly brought their friends home to admire his ``ghastly generalized manifolds.'' They understood almost none of the mathematics, but they immediately understood Bob's enthusiasm. Of all the mathematicians who visited our home, he is the one they remember most vividly.

Working with Bob never felt routine. A discussion that began with one problem might end several hours later with several new questions, two possible counterexamples, and the outline of another paper. Our collaboration led to work on shrinking criteria and on general-position properties characterizing 3-manifolds. Some of the questions that occupied us remain open, especially in dimension four. Yet what I remember most vividly is not the finished papers, but the creative momentum that produced them. Bob never seemed interested in proving that he was the cleverest person in the room. He wanted everyone in the room to enjoy mathematics as much as he did, and he had a rare ability to make others believe that they, too, could contribute.

Bob was also tremendous fun to travel with. In 1989 we went by train to the Colloquium on Topology in P\'ecs, Hungary, organized by the J\'anos Bolyai Mathematical Society. The meeting brought together many topologists from countries then still behind the Iron Curtain and gave Bob an opportunity to meet many Soviet topologists in person. After crossing into Hungary, we transferred to an extraordinarily slow local train. It seemed to stop at every village so that passengers carrying large bags could get on and off. After yet another long halt, Bob observed dryly that we would probably reach P\'ecs sooner if we had gotten off at the border and just walked. For years afterward, that journey remained one of our favorite shared memories.

A story about our Soviet publishing venture from that period also made us laugh for many years afterward. In 1987 I presented our joint work at a very large international topology conference in the Soviet Union, which took place in Baku, Azerbaijan. The translation of our paper was published in the Russian Steklov Mathematical Institute Proceedings in 1992, and the English version was to be published by the American Mathematical Society one year later. We were never sent any galley proofs for the Russian translation which unfortunately, contained some errors, most important was in the title -  instead of {\it 1-demensional decompositions} (i.e. {\it of embedding dimension 1}) they put
{\it 1-dimensional decompositions} (i.e. {\it of covering dimension 1}). Therefore we had to spend a lot of effort to make sure that the AMS would use our original English manuscript instead of by default, translating the Russian version. Also, since in the Soviet Union authors were entitled to an honorarium, Bob suggested many years later that some roubles might still be waiting for us in Moscow and that the accrued interest could now make it worth claiming them.

\begin{figure}[ht]
\centering
\includegraphics[width=\columnwidth]{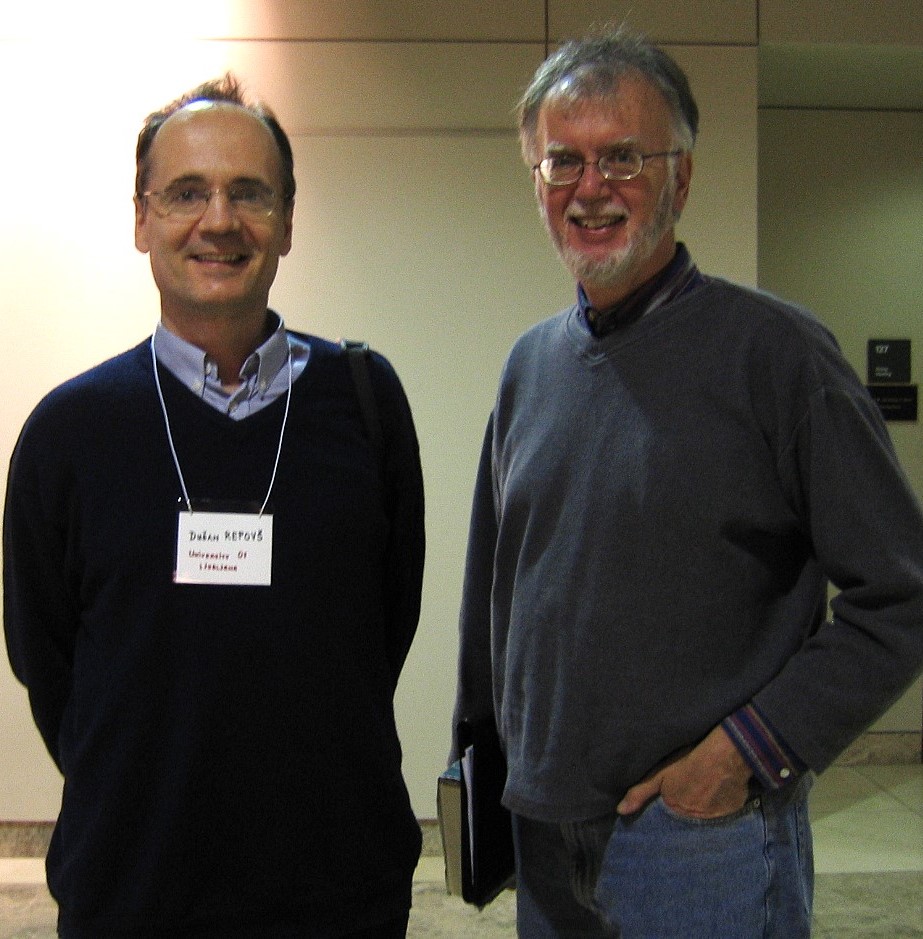}
\caption{Du\v{s}an Repov\v{s} with Robert J. Daverman at the Barrett Lectures, University of Tennessee, Knoxville, 2006.}
\end{figure}

Bob's humor could also turn gently toward his own subject. Although decomposition theory was one of his lifelong passions, he once remarked that decomposition theorists should not work too hard to persuade others of its importance: everyone in the field already believed in it \emph{a priori}, and usually far too strongly. That combination of devotion and self-mockery was characteristic. He took mathematics seriously, but never himself too seriously.

I last visited Bob in Knoxville in 2006, for the Barrett Lectures. With him, mathematical conversations did not really end; they merely paused and resumed years later. His last email to me concerned his 2019 paper with Thomas Thickstun on degree-one monotone self-maps of the Pontryagin surface. Their paper answered affirmatively a question about these exotic spaces that I had raised several times in the 1990s. I was delighted by their success, and equally delighted that a conversation begun decades earlier had found another continuation.

My wife Barbara always delighted in welcoming visiting mathematicians into our home. She was an outstanding cook, following a family tradition inherited from her grandmother, who had spent her life as a professional chef, and she took particular pleasure in preparing a feast for all our guests. Over the years Barbara met many of my mathematical friends, but she always looked forward especially to Bob's visits. She would say that he was the most amusing of them all and then, after a thoughtful pause, add that he was also ``the easiest to talk to.'' I have never found a better description of him.

\end{document}